\documentstyle[11pt]{article}
\newlength{\spacing}
\newcommand{\doublespace}{\setlength{\baselineskip}{1.0\spacing}}
\newtheorem{thm}{Theorem}[section]

\newtheorem{lem}[thm]{Lemma}
\newtheorem{prop}[thm]{Proposition}
\def\lam{\lambda }

\def\lab{\label }

\def\rar{\to}

\def\inft{\infty}

\def\del{\delta}

\def\ep{\epsilon}

\def\lam{\lambda}
\def\vep{\varepsilon}

\def\today{\ifcase\month\or
  January\or February\or March\or April\or May\or June\or
  July\or August\or September\or October\or November\or December\fi
  \space\number\day, \number\year}

\begin{document}
\begin{titlepage}
\begin{center}
{\bf Vertex Degree of Random Geometric Graph on Exponentially Distributed Points} \\
\vspace{0.20in} by \\
\vspace{0.2in} {Bhupendra Gupta \footnote{Corresponding Author.
email:gupta.bhupendra@gmail.com, bhupen@iitk.ac.in}}\\
Department of Mathematics, Indian Institute of Technology,
Kanpur 208016, India \\
\vspace{0.1in}
\end{center}
\vspace{0.2in}
%
%\begin{center} {\bf Preprint} \end{center}
\sloppy
\begin{center} {\bf Abstract} \end{center}

\begin{center} \parbox{4.8in}
{Let $X_1,X_2,\ldots$ be an infinite sequence of i.i.d. random
vectors distributed exponentially with parameter $\lam .$ For each
$y$  and $n\geq 1,$ form a graph $G_n(y)$ with vertex set $V_n =
\{X_1,\ldots,X_n\},$ two vertices are connected if and only if edge
distance between them is greater then $y$, i.e, $\|X_i-X_j\| \leq
y.$ Almost-sure asymptotic rates of convergence/divergence are
obtained for the minimum and maximum vertex degree of the random
geometric graph, as the number of vertices becomes large $n,$ and
the edge distance varies with the number of vertices.
%Let $n$ points be placed independently in $d-$dimensional space
%according to the density $f(x) = C_d\lam^d e^{-\lam
%\|x\|},\;\lam>0,\;x \in {\Re}^d,\;d \geq 2.$ Let $d_n$ be the
%longest edge length for the nearest neighbor graph on these points.
%We show the
%strong law result, $\lam (\log \log n)^{-1/2}d_n \rar d$ as $n \rar \infty$ a.s.
} \\
\vspace{0.4in}
\end{center}

%\rule[1mm]{4in}{.4mm} \\
\vspace{0.5in}
{\sl AMS 1991 subject classifications}: \\
\hspace*{0.5in} Primary:   60D05, 60G70\\
\hspace*{0.5in} Secondary:  05C05, 90C27
\end{titlepage}
\doublespace
\section{Introduction\lab{s1}}
The development in this paper parallels that of Appel and Russo
\cite{Appel,Appel1}. Let $d \geq 1$ be an integer. %Let
%the matric on $\mathbb{R}^d$ be given by an arbitrary specified norm
%$\|.\|.$
Some strong law results for maximum and minimum vertex degree for
the uniformly distributed points on the d-dimensional unit cube,
using the ${\it l}_{\infty}$ norm, are given by Appel and Russo
(\cite{Appel}, \cite{Appel1}). Also Penrose gave some convergence in
probability results in his book \cite{Penrose}. A strong law result
for the maximum vertex degree of a graph whose vertices from a
density with compact support is given by Penrose \cite{Penrose}. A
considerable amount of work has been done in this field, but most of
the the work is for densities with compact support and in most of
the cases when the density is uniform. It is more relevant to
consider cases such as in which distribution of vertices have
unbounded support, like the {\it weak law} result for the {\it
nearest neighbor distance} for the normally distributed vertices
given by M.Penrose\cite{Penrose1}. Gupta and Iyer \cite{Gupta} give
{\it strong law} result for the {\it largest nearest neighbor
distance} $d_n$ for the normally distributed vertices. Also, Gupta
and Iyer \cite{Gupta1} prove that when the tail of the density
decays like an exponential or slower, $d_n$ diverges, whereas for
super exponential decay of the tail, $d_n \rar 0,$ a.s. as $n \rar
\infty.$ Properties of the one dimensional exponential random
geometric graphs have been studied in Iyer, Manjunath and Gupta
(2005). In the one dimensional exponential case, the spacings
between the ordered nodes are independent and exponentially
distributed. This allows for explicit computations of many
characteristics for the graph and both strong and weak law results
can be established. For a detailed description of Random Geometric
Graphs, their properties and applications, we
refer the reader to Penrose \cite{Penrose} and references therein.\\

Unlike the normal distribution, there is no unique natural extension
of the one dimensional exponential distribution. So we consider the
following criteria, let $X_1,X_2,X_3,\ldots$ be independent
d-dimensional random vectors and distributed according to the
exponential density along each axis. Let the vertex set ${\cal{X}}_n
= \{X_1,X_2,\ldots X_n\}, \: n = 1, 2, 3, \ldots$ be  the point
process. This work concerns with the maximum and minimum vertex
degree of the random geometric graph $G_n(y)$. We begin our
investigation of the structure of the random graph $G_n(y)$ for
various values of $y$ considering the edge set $E_n(y) =
\{\{X_i,X_j\}:\|X_i-X_j\| \leq y, 1 \leq i \neq j \leq n\}.$ We
obtain strong law results as did Appel
and Russo \cite{Appel,Appel1}.\\
\section{Supporting Results\lab{s2}}
\begin{lem}
$X$ follow a binomial distribution with parameters $n$ and $p.$ Then
\begin{equation}
P[X \geq k] \leq \left(\frac{np}{k}\right)^k \exp(k-np), \qquad
k\geq np,\label{e1}
\end{equation}
and
\begin{equation}
P[X \leq k] \leq \left(\frac{np}{k}\right)^k \exp(k-np), \qquad 0
\leq k\leq np.\label{e2}
\end{equation}
 \label{lemma1}
\end{lem}
For convenience, we may express the bounds in (\ref{e1}) and
(\ref{e2}) as
\[\exp\left(npH\left(\frac{np}{k}\right)\right),\]
where
\[H(t) = \frac{1}{t}\log\:t + \frac{1}{t} -1,\quad 0<t<\inft,\]
and $H(\inft) = -1.$ Note that $H(t)<0$ for all $0<t \neq 1 \leq
\infty ;\: H$ is increasing on $(0,1)$ and decreasing on
$(1,\infty).$\\

{\bf Proof.} See Lemma 1.2 on page 25, Penrose \cite{Penrose}.\\

Let $X_1,X_2,\ldots$ be an infinite sequence of i.i.d. random
vectors distributed exponentially with parameter $\lam .$ For each
$y$  and $n\geq 1,$ form a graph $G_n(y)$ with vertex set $V_n =
\{X_1,\ldots,X_n\},$ in which two vertices are connected if and only
if distance between them is less then $y$, i.e, $\|X_i-X_j\| \leq
y.$ We refer to the parameter $y$ as the edge distance.\\
 We begin our investigation of the structure of the random graph
 $G_n(y)$ for various values of $y$ considering the edge set
$E_n(y) = \{\{X_i,X_j\}:\|X_i-X_j\| \leq y, 1 \leq i \neq j \leq
n\}.$ Let $\vep_n(y) = card[E_n(y)],$ the number of the edges in the
random graph $G_n(y).$ Then
\begin{equation}
p(y) = P[\|X_i-X_j\|\leq y] = (1- e^{-\lam y})^d.
\label{connectivity of two node}
\end{equation}
We note that $p(y) \approx (\lam y)^d,$ as $y\downarrow 0,$ the
volume of an $l_{\inft}$-ball of radius y.\\

For fixed $y,$ the number of the edges $\vep_n(y)$ in the random
graph is ${n\choose 2} U,$ where $U$ is a U-statistics with
symmetric kernel $h(u,v) = 1_{\{\|u-v\|\leq y\}}$ and expected value
is
\begin{equation}
E[\vep_n(y)] = {n\choose 2}p(y).
\end{equation}
Standard theory (Serfling \cite{Serfling}, Chapter 5) provides a
strong law of large numbers:
\begin{equation}
\frac{\vep_n(y_n)}{{n\choose 2}} \rar p(y_n)\quad a.s., \qquad as\:
n\rar \infty.\label{e3}
\end{equation}
This convergence can easily proved to be uniform by the standard
proof of the Glivenko-Cantelli Theorem (See Chow and Teicher
\cite{chow}, Section 8.2).
\begin{prop}
\begin{equation}
 \sup_{0\leq y \leq 1}\left|\frac{\vep_n(y)}{{n\choose 2}}-
p(y)\right| \rar 0 \qquad a.s. \qquad as \: n\rar \infty.
\end{equation}
\label{pro1}
\end{prop}
We can prove the following rate of convergence in the strong law
(\ref{e3}) along a particular sequence of edge distances.
\begin{prop} If $\{y_n\}$ be any sequence of edge distances
that satisfies
\begin{equation}
\frac{ny_n^d}{\log{n}} \rar c \in (0,\inft], \label{seq of edge
distance}
\end{equation}
for large $n.$ Then,
\begin{equation}
\lim_{n\rar \inft}\left|{\frac{\vep_n(y_n)}{{n\choose 2}}}-
p(y_n)\right|= 0 \qquad a.s.
\end{equation}\label{proposiotion for edge distances}
\end{prop}
{\bf Proof.} See Proposition 2, Appel and Russo, \cite{Appel}.\\

The following lemma gives the region of interest when the points are
independently exponentially distributed along the axes.
For sufficiently large $n,$ define
\begin{equation}
R_n = \frac{\log\:n}{\lam d}.
\end{equation}
Let $U_n$ be the event ${\cal{X}}_n \subset B(0,R_n)$ i.e. all the
points of ${\cal{X}}_n$ lie in $B(0,R_n).$
\begin{lem}
For sufficiently large $n,$
\[P[U_n^c \:\:i.o.] = 0 \qquad  a.s.\]
\label{lemma2}
\end{lem}
{\bf Proof.} \[P[U_n^c] = P[\cup_{i=1}^{n}(X_i > R_n)].\]
Let $n_k$ be the subsequence $k^a$, with $a>0,$ and consider
\begin{eqnarray*}
P[\cup_{n=n_k}^{n_{k+1}}U_n^c] & \leq & P[\cup_{i = 1}^{n_{k+1}}(X_i
> R_{n_k})]\\
& \leq & \sum_{i = 1}^{n_{k+1}}P[(X_i > R_{n_k})]\\
& = & n_{k+1}(\exp(-\lam R_{n_k}))^d,
\end{eqnarray*}
take $R_n = \frac{1+\ep}{\lam d}\log\:n$
\begin{eqnarray}
P[\cup_{n=n_k}^{n_{k+1}}U_n^c] & \leq &
\frac{n_{k+1}}{(n_k)^{1+\ep}}\nonumber\\
& \approx & \frac{1}{k^a\ep}.\label{R_n}
\end{eqnarray}
We can always choose $a$ sufficiently large such that (\ref{R_n}) is
summable. Now by the Borel-Cantelli Lemma, we get all the vertex of
${\cal{X}}_n$ are lie in $B(0,R_n)$ a.s. eventually for $R_n =
\frac{1}{\lam d}\log\:n.$
\section{Minimum Vertex Degree\lab{s3}}
For $1 \leq i\leq n,$ let
\begin{equation}
deg X_{n,i} = \sum_{1\leq j \neq i \leq n}1_{\{\|X_i - X_j\| \leq
x\}},
\end{equation}
be the degree of the vertex $X_i$ in the random graph $G_n(y)$ and
let
\begin{equation}
\delta_n(y_n) = \min \{deg X_{n,1}(y_n),\ldots,deg X_{n,n}(y_n)\},
\end{equation}
be the minimum vertex degree of $G_n(y_n).$ Note that $\del_n$ is
non-decreasing in $y_n$ for fixed $n.$
\begin{thm}
Let $\del_n$ be the minimum vertex degree of the $G_n$ defined on
the collection ${\cal{X}}_n$ of $n$ points distributed independently
and identically according to the exponential density $f(\cdot).$
Then,
\begin{equation}
\limsup_{n \rar \infty}\frac{\del_n(y_n)}{ny_n^d} \leq (\lam)^d,
\qquad a.s.,
\end{equation}
and
\begin{equation}
\liminf_{n \rar \inft}{\frac{\delta_n(y_n)}{n y_n^d}} \geq
a(c)\lam^d, \qquad a.s.,
\end{equation}
where $a(c)$ is the  root in $(0,1)$ of
\begin{equation}
a\log a -a +1 = {\frac{1}{\lam^d c}},
\end{equation}
with $a(\inft) = 1.$
\label{theorem1}
\end{thm}
{\bf Proof.} Pick $u,t$ and $\ep>0,$ such that
\[\ep +u = t.\]
From Lemma \ref{lemma2}, ${\cal{X}}_n \subset B(0,R_n)$ a.s. for all
large enough $n.$ For fixed $t>0,$ define $S_n(t) = B(R_n,y_n(u)),$
and define $y_n(\cdot)$ as follows
\begin{equation}
\lam y_n(t) = \left(\frac{t \log\:n}{n}\right)^{1/d}. \label{x_n}
\end{equation}
Define $B(R_n,y_n(u)) = B(R_n e,y_n(u)),$ where $`e'$ is the
d-dimensional unit vector $(1,0,\ldots,0).$ Define an event $W_n$
such that there is a vertex $X_1$ of ${\cal{X}}_n$ lie in
$B(R_n,y_n(\ep))$ and at least one another vertex of ${\cal{X}}_n$
is lie in $B(R_n,y_n(u))$ i.e., in $y_n(u)$ neighborhood of $X_1$.

\[W_n =  \cup_{j=1}^{n-1}\{X_j \in [0,R_n]^d,\: x_n \in
B(R_n, y_n(\ep))\: : \|x_n -X_ j\| \leq y_n(u)\}.\]
Thus,
\begin{eqnarray*}
P[W_n] & = & P[\cup_{j=1}^{n-1}\{X_j \in [0,R_n]^d,\: x_n
\in B(R_n, y_n(\ep))\: : \|x_n -X_ j\| \leq y_n(u)\}]\\
& \leq & \sum_{j=1}^{n-1}P[X_j \in [0,R_n]^d,\: x_n \in B(R_n, y_n(\ep))
\: : \|x_n -X_ j\| \leq y_n(u)]\\
& = & (n-1)P[X_j \in [0,R_n]^d,\: x_n \in B(R_n, y_n(\ep))\: :
\|x_n -X_ j\| \leq y_n(u)]\\
& \sim & (n-1)(e^{-\lam(R_n-y_n(u))} - e^{-\lam(R_n+y_n(u))})^d\\
& = & (n-1) e^{-\lam d R_n}(e^{\lam y_n(u)}-e^{-\lam y_n(u)})^d\\
& = & (n-1) e^{-\lam d R_n}(1+\lam y_n(u) - 1 +\lam y_n(u))^d\\
& \sim & \frac{(n-1)(2\lam)^d y_n^d(u)}{n} \sim (2\lam)^d y_n^d(u)\\
& := & q_n.
\end{eqnarray*}
Let $Y_n(t) = \sum_{i=1}^{n}I_{\{X_i \in S_n(t)\}}$ denote the
number of the vertices in sub-cube $S_n(t).$ Note that $Y_n(t)$ is
distributed according to $Bi(n,q_n),$ where $q_n$ is as defined above.\\

Let $b > u$ be given. Taking a subsequence such the $n_k = k^a,$
where $a>0.$
\begin{equation}
P[\cup_{n = n_k}^{n_{k+1}}\{Y_n > b \log\:n\}]  \leq  P[X_{n_k}>
b\log\:n],
\end{equation}
where $X_{n_k} = \sum_{i=1}^{n_{k+1}}I_{\{X_i \in S_{n_k}(t)\}}.$
\begin{eqnarray}
P[\cup_{n = n_k}^{n_{k+1}}\{Y_n > b \log\:n\}] & \leq & (k+1)^a
\exp\left((k+1)^a q_{n_k}H\left(\frac{(k+1)^a
q_{n_k}}{b\log\:n_k}\right)\right)\nonumber\\
& \leq & (k+1)^a \exp\left((k+1)^a (2\lam)^d
y_{n_k}^d(u)H\left(\frac{(k+1)^a(2\lam)^d
y_{n_k}^d(u)}{b\log\:n_k}\right)\right)\nonumber\\
& \sim & (k+1)^a \exp\left((k+1)^a \frac{2^d ua\log\:k}{k^a}
H\left(\frac{(k+1)^a 2^d ua \log\:k}{ba\log\:k\: k^a}\right)\right)\nonumber\\
& \sim & (k+1)^a \exp\left(2^d ua\log\:k H\left(\frac{2^d
u}{b}\right)\right),
\end{eqnarray}
since $b > u, \: H(2^d u/b)\leq 0.$ Hence
\begin{eqnarray}
P[\cup_{n = n_k}^{n_{k+1}}\{Y_n > b \log\:n\}] & \leq & (k+1)^a
\exp(-Cau\log\:k)\nonumber\\
& = & \frac{(k+1)^a}{k^{2^d Cau}} \sim \frac{1}{m^{(2^dCu-1)a}},
\end{eqnarray}
where $C = - H(2^du/b).$ If $2^dCu>1,$ then $H(2^du/b)<-1.$ Then for
sufficiently large $a$ the above probability is summable. Then by
the the Borel-Cantelli Lemma, we get
\begin{equation}
Y_n(t) \leq b\log\:n \qquad a.s.,\: \mbox{eventually, for any $b>
u.$}\label{y_n}
\end{equation}
This implies the $Y_n(t),$ number of vertices inside the ball
$B(R_n,y_n(u))$ is at most $b\log\:n$ a.s.\\
Finally, Let $\ep >0,$ be given and let
\[V_n = \cap_{i=1}^{n}\{X_i \notin S_n(\ep)\}.\]
The monotonicity of edge length implies that
\begin{eqnarray*}
P[\cup_{n = n_k}^{n_{k+1}}V_n] & \leq &  P[\cap_{i=1}^{n^k}\{X_i
\notin S_n(\ep)\}]\\
& = & \left(\exp\left(-\frac{\ep\log n_{k+1}}{(2\lam)^dn_{k+1}}
\right)\right)^{n_k}\\
& = & \exp\left(-\frac{C_1k^aa\log(k+1)}{(k+1)^a}\right)\\
& \sim & \exp(-C_1a\log(k+1))\\
& = & \frac{1}{(k+1)^{aC_1}},
\end{eqnarray*}
we can always choose $a$ sufficiently large, such that
\[\sum_{k=0}^{\infty} P[\cup_{n = n_k}^{n_{k+1}}V_n] < \infty.\]
Therefore by the Borel-Cantelli Lemma, for large enough $n,$ there
is atleast one vertex $X_i$ of ${\cal{X}}_n$ in the sub-cube
$S_n(\ep),$ a.s. for some $i= 1,2,\ldots,n.$ For fixed $u>0$ any
such vertex  has, a.s., at most $u\log\:n$ neighbors at edge
distance $y_n(u)$ eventually, according to (\ref{y_n}). Thus,
\begin{equation}
\limsup_{n \rar \infty}\frac{\del_n(y_n)}{ny_n^d(u)} \leq (2\lam)^d,
\qquad a.s.
\end{equation}
\\

%
%\begin{thm}
% Let $\{y_n\}$ be any sequence of edge distances which satisfies
% connectivity regime define as in equation(\ref{seq of edge distance}). Then
% %
% \begin{equation}
%\liminf_{n \rar \inft}{\frac{\delta_n(y_n)}{n y_n^d}} \geq
%a(c)\lam^d, \qquad a.s.,
%\end{equation}
%%
%where $a(c)$ is the  root in $(0,1)$ of
%\begin{equation}
%a\log a -a +1 = {\frac{1}{\lam^d c}}
%\end{equation}
%with $a(\inft) = 1.$
%\end{thm}
To established the limit infimum of minimum vertex degree. Consider
\begin{eqnarray*}
P[\del_n(y_n)\leq K] & = &
P[\min\{degX_{n,1}(y_n),\ldots,degX_{n,n}(y_n)\}\leq K]\\
& = & P[\cup_{i=1}^{n}(degX_{n,i}(y_n)\leq K)]\\
& \leq & \sum_{i=1}^{n}P[degX_{n,i}(y_n)\leq K)]\\
& \leq & n P[degX_{n,i}(y_n)\leq K)].
\end{eqnarray*}
Since $degX_{i} \sim bi(n,p(y_n)),$ where $p(y_n)$ is define as in
the equation(\ref{connectivity of two node}). Then by the
Chernoff's bound
\begin{equation}
P[\del_n(y_n) \leq K] \leq n\exp \left(n(\lam y_n)^d
H\left(\frac{n(\lam y_n)^d}{K}\right)\right).
\end{equation}
Let $a \in (0,c)$ be given and let
\[u_n^d ={\frac{(c-a)\log n}{n}} \qquad n\geq 1,\]
 and
\[v_n^d ={\frac{(c+a)\log n}{n}}\qquad n\geq 1.\]
The edge distances $\{u_n\}$ and $\{v_n\}$ are both decreasing in
$n$ with $u_n\leq y_n\leq v_n$ for all $n$ large enough. Then
\begin{eqnarray*}
P\left[\cup_{n={n_k}}^{n_{k+1}}\left(\frac{\del_n(v_n)}{nu_n^d} \leq
\ep\right)\right] & = &
P\left[\frac{\del_{n_{k+1}}(v_{n_k})}{n_{k+1}u_{n_{k+1}}^{d}} \leq
\ep \right]\\
& \leq &
n_{k+1}\exp\left(n_{k+1}p(v_{n_k})H\left(\frac{n_{k+1}p(v_{n_k})}{\ep
n_k u_{n_{k+1}}^d}\right)\right)\\
& = & n_{k+1}\exp\left(\lam^d (c+a)\log n_k
H\left(\frac{\lam^d(c+a)}{\ep(c-a)}\right)\right),
\end{eqnarray*}
as $a \rar 0,$
\begin{eqnarray}
P\left[\cup_{n = n_k}^{n_{k+1}}\left(\frac{\del_n(v_n)}{nu_n^d} \leq
\ep \right)\right] & \leq & n_{k+1}\exp\left(\lam^d c \log
n_k H\left(\frac{\lam^d}{\ep}\right)\right)\nonumber\\
& = & n_{k+1}(n_k)^{\eta},
\end{eqnarray}
where $\eta = \lam^d cH\left(\frac{\lam^d}{\ep}\right).$ The above
expression is summable if $\eta < -1 \: i.e., \: \ep \leq
a(c)\lam^d,$ where $a(c)$ is decreasing in $(0,1)$ and defined as
in the statement. Then by the Borel-Cantelli Lemma, we get the
required result.
\section{Maximum Vertex Degree\lab{s4}}
For $1 \leq i\leq n,$ let
\begin{equation}
deg X_{n,i} = \sum_{1\leq j \neq i \leq n}1_{\{\|X_i - X_j\| \leq
x\}},
\end{equation}
be the degree of the vertex $X_i$ in the random graph $G_n(y)$ and
let
\begin{equation}
\Delta_n(y) = \max \{deg X_{n,1}(y_n),\ldots,deg X_{n,n}(y_n)\},
\end{equation}
be the maximum vertex degree of $G_n(x).$ One of our main aim is,
for any particular vertex, find out the rate at which
$\Delta_n(y_n)$ diverges, for all $\{y_n\}$ throughout the
connectivity regime define as in (\ref{seq of edge distance}).
\begin{thm}
Let $\{y_n\}$ be a sequence of the edge distances which satisfies
connectivity regime define as in (\ref{seq of edge distance}). Then
\begin{equation}
\liminf_{n \rar \inft}{\frac{\Delta_n(y_n)}{n y_n^d}} \geq \lam^d,
\qquad a.s.,
\end{equation}
\begin{equation}
\limsup_{n \rar \inft}{\frac{\Delta_n(y_n)}{n y_n^d}} \leq
a(c)\lam^d, \qquad a.s.,
\end{equation}
where $a(c)$ is the  root in $[1,\inft)$ of
\begin{equation}
a\log a -a +1 = {\frac{1}{\lam^d c}},
\end{equation}
with $a(\inft) = 1.$
\end{thm}
Proof. We know that $2\vep_n(y_n) \leq n \Delta_n(y_n),$ and so
\begin{eqnarray*}
\frac{\Delta_n(y_n)}{(n-1)y_n^d} & \geq &
\frac{\vep_n(y_n)}{{n\choose 2} y_n^d}\\
& = & \frac{p(y_n)}{y_n^d} +
\frac{1}{y_n^d}\left(\frac{\vep_n(y_n)}{{n \choose 2}} -
p(y_n)\right)\\
& = & \lam^d,
\end{eqnarray*}
by using the proposition(\ref{proposiotion for edge distances}). The
above expression implies the
\begin{equation}
\liminf_{n \rar \inft} \frac{\Delta_n(y_n)}{n y_n^d} \geq \lam^d.
\end{equation}
\\
To prove second part, consider
\begin{eqnarray}
P[\Delta_n(y_n) \geq K] & = & P[\max\{deg X_{n,1}(y_n),\ldots,deg
X_{n,n}(y_n)\} \geq K]\nonumber\\
& = & P[\cup_{i=1}^{n}(deg X_{n,i}(y_n)\geq K)]\nonumber\\
& \leq & \sum_{i=1}^{n}P[deg X_{n,i}(y_n)\geq K]\nonumber\\
& = & n P[deg X_{n,i}(y_n)\geq K].
\end{eqnarray}
Since $deg X_i \sim bi(n,p(y_n)),$ where $p(y_n)$ is define as in
equ.(\ref{connectivity of two node}). Then by the Dudley's
inequality
\begin{equation}
P[\Delta_n(y_n)\geq K] \leq n
\exp\left((n/2)p(y_n)H\left({\frac{np(y_n)}{2K}}\right)\right).
\end{equation}
Let $a \in (0,c)$ be given and let
\[u_n^d ={\frac{(c-a)\log n}{n}} \qquad n\geq 1,\]
 and
\[v_n^d ={\frac{(c+a)\log n}{n}}\qquad n\geq 1.\]
The edge distances $\{u_n\}$ and $\{v_n\}$ are both decreasing in
$n$ with $u_n\leq y_n\leq v_n$ for all $n$ large enough. Then
\begin{eqnarray}
P\left[\cup_{n_k}^{n=n_{k+1}}\left(\frac{\Delta_n(v_n)}{nu_n^d} \geq
\ep \right)\right] & = & P\left[{\frac{\Delta_{n_k+1}(v_{n_k})}{n_k
u^d_{n_{k+1}}}} \geq
\ep \right]\nonumber\\
& \leq & n_{k+1} \exp \left(n_{k+1}p(v_{n_k})
H\left({\frac{n_{k+1}p(v_{n_k})}{\ep n_k
u^d_{n_{k+1}}}}\right)\right)\nonumber\\
& = & n_{k+1} \exp \left(\lam^d(c+a) \log n_k H\left(
{\frac{\lam^d(c+a)}{\ep(c-a)}}\right)\right),
\end{eqnarray}
as $a \rar 0,$
\begin{eqnarray}
P\left[\cup_{n_k}^{n=n_{k+1}}\left(\frac{\Delta_n(v_n)}{nu_n^d} \geq
\ep \right)\right] & \leq & n_{k+1} \exp \left(\lam^dc \log n_k
H\left(
{\frac{\lam^d}{\ep}}\right)\right)\nonumber\\
& = & n_{k+1} (n_k)^\eta,
\end{eqnarray}
where $\eta = \lam^d c H({\frac{\lam^d}{\ep}})$, the above
expression will be summable if $\eta <-1$ i.e., $\ep \leq
a(c)\lam^d,$ where $a(c)$ is increasing in $[1,\infty)$ and
defined as in the statement. Then by the Borel-Cantelli Lemma, we
get
\begin{equation}
P\left[\cup_{n_k}^{n=n_{k+1}}\left(\frac{\Delta_n(v_n)}{nu_n^d} \geq
\ep \right)\qquad i.o.\right] = 0.
\end{equation}
This implies that,
\begin{equation}
\limsup_{n \rar \inft}\frac{\Delta_n(y_n)}{ny_n^d} \leq a(c)\lam^d,
\qquad a.s.
\end{equation}
\begin{thm}
 Let $\{y_n\}$ be any sequence of edge distances and let $S = \sum_n n
 y_n^d.$\\
 1. If $S < \inft$ and $y_n$ is non increasing then $P[\Delta_n(y_n) \geq 1, i.o.] =
 0.$\\
 2. If $S = \inft$ then $P[\Delta_n(y_n) \geq 1, i.o.] =
 1.$
\end{thm}
Proof. We note that the expected degree of a given vertex in the
random graph $G_n(y)$ is $(\lam y_n)^d.$\\
First we consider part 1. Assume $S < \inft$ and that $y_n$ is non
increasing.\\
Let \[R_n(y) = \cup_{j=1}^{n-1}\{x \in [0,\inft]^d : \|U_j -x\| \leq
y\}\] the union of the balls of radius $y$ centered at first
$(n-1)$ vertices.\\
Let $E_n = \{U_n \in R_n(x_n)\}.$ Then by the Boole's inequality
\begin{eqnarray*}
P[E_n] & = & P[U_n \in R_n(y_n)]\\
& = & P[\cup_{j=1}^{n-1}\{U_n \in [0,\inft]^d : \|U_j -U_n\|
\leq y_n\}]\\
& \leq & \sum_{j=1}^{n-1}P[U_n \in [0,\inft]^d : \|U_j -U_n\|
\leq y_n]\\
& = & (n-1)(1-e^{-\lam y_n})^d\\
& \approx & \lam^d (n-1) y^d_n.
\end{eqnarray*}
This probability summable, i.e. $\sum_{n=0}^{\inft}P[E_n] <\inft.$
Hence by the Borel-Cantelli Lemma,
\begin{equation}
P[E_n,\quad \mbox{i.o.}] = 0.
\end{equation}
For each $\omega \in \{E_n,\quad \mbox{i.o.}\}, \exists $ a positive
integer $N = N(\omega)$ such that $\|U_j -U_n\| > y_n$ for each $j =
1,\ldots, n-1,$ Whenever $n \geq N.$ Since $y_n \rar 0,$ we may
choose $N_1 >N$ such that $\|U_j -U_n\| > y_{N_1}$ for each $i,j =
1,\ldots, n-1.$ Given $n > N_1,$ we observe that for $1\leq j< k
\leq n,\:k > N$ implies that $\|U_j -U_n\| > y_k \geq y_n$ while $k
> N$ implies that $\|U_j -U_n\| > y_{N_1} > y_n.$
Hence $G_n(y_n)$ has no edge.\\

Now, we consider part 2,\\
Assume that $S = \inft$ and
\[R_n(y) = \cup_{j=1}^{n-1}\{x \in [0,\inft]^d : \|U_j -x\|
\leq y\}.\]\\
Let $ D_n = \{R_n(y_n/2) \mbox{is a disjoint union}\}$ and let $D =
\{R_n(y_n/2) \mbox{is disjoint ult.}\} = \cup_M \cap_{m \geq
M}D_m.$\\
Let $\{{\cal{F}}_n\},$ be a family of sigma algebras ${\cal{F}}_n=
\sigma(U_1,\ldots, U_n)$. Note that $\omega \in D_n.$
\[P[U_n \in R_n(y_n/2)\mid {\cal{F}}_{n-1}](\omega) =
 (n-1)(\lam y_n/2)^d \approx (\lam /2)^d n y_n^d,\]
for large $n.$ Thus on $D$
\[\sum_n P[U_n \in R_n(y_n/2)\mid {\cal{F}}_{n-1}] = \inft, \]
since $S = \inft.$\\
Hence by the conditional Borel-Cantelli Lemma
\[P[U_n \in R_n(y_n/2)\mid {\cal{F}}_{n-1},\qquad \mbox{i.o.}] = 1.\]
This implies that
\[P[U_n \in R_n(y_n/2),\qquad \mbox{i.o.}] = 1.\]
Hence
\[P[\Delta_n(y_n) \geq 1, \qquad \mbox{i.o.}] = 1.\]
\end{document}